\documentclass[amsmath,amssymb,aps,prb,twocolumn,superscriptaddress]{revtex4-2}
\usepackage{graphicx}
\usepackage{bm}

\begin{document}


\title{Geometry-informed dynamic mode decomposition in origami dynamics}


\author{Shuaifeng Li}
\affiliation{Department of Aeronautics and Astronautics, University of Washington, Seattle, WA, USA}
\author{Yasuhiro Miyazawa}
\affiliation{Department of Aeronautics and Astronautics, University of Washington, Seattle, WA, USA}
\author{Koshiro Yamaguchi}
\affiliation{Department of Aeronautics and Astronautics, University of Washington, Seattle, WA, USA}
\author{Panayotis G.Kevrekidis}
\affiliation{Department of Mathematics and Statistics, University of Massachusetts, Amherst, MA, USA}
\author{Jinkyu Yang}
\affiliation{Department of Aeronautics and Astronautics, University of Washington, Seattle, WA, USA}
\affiliation{Department of Mechanical Engineering, Seoul National University, Seoul, Republic of Korea}


\begin{abstract}
Origami structures often serve as the building block of mechanical systems due to their rich static and dynamic behaviors. Experimental observation and theoretical modeling of origami dynamics have been reported extensively, whereas the data-driven modeling of origami dynamics is still challenging due to the intrinsic nonlinearity of the system. In this study, we show how the dynamic mode decomposition (DMD) method can be enhanced by integrating geometry information of the origami structure to model origami dynamics in an efficient and accurate manner. In particular, an improved version of DMD with control, that we term geometry-informed dynamic mode decomposition~(giDMD), is developed and evaluated on the origami chain and dual Kresling origami structure to reveal the efficacy and interpretability. We show that giDMD can accurately predict the dynamics of an origami chain across frequencies, where the topological boundary state can be identified by the characteristics of giDMD. Moreover, the periodic intrawell motion can be accurately predicted in the dual origami structure. The type of dynamics in the dual origami structure can also be identified. The model learned by the giDMD also reveals the influential geometrical parameters in the origami dynamics, indicating the interpretability of this method. The accurate prediction of chaotic dynamics remains a challenge for the method. Nevertheless, we expect that the proposed giDMD approach will be helpful towards the prediction and identification of dynamics in complex origami structures, while paving the way to the application to a wider variety of lightweight and deployable structures.
\end{abstract}


\maketitle

\section{Introduction\label{introduction}}
Origami, as an ancient handcrafted paper folding art, captivates not only artists and mathematicians with its exquisite design principles, but also engineers with its enormous (and increasingly leveraged in recent years) potential in engineering applications. By introducing the principles of creasing and folding to the flat materials, origami structures can be formed, which lead to stiffness enhancement~\cite{filipov2015origami}, negative Poisson’s ratio~\cite{lv2014origami,yasuda2015reentrant} and multistability in the architecture~\cite{jianguo2015bistable,hanna2014waterbomb,yasuda2017origami}. These developments also inspire further applications in the robotics~\cite{miyashita2015untethered,zirbel2013accommodating}, medical equipment~\cite{edmondson2013oriceps,nelson2016curved}, and mechanical metamaterials with unprecedented mechanical properties~\cite{deng2022tunable,lyu2021origami}. Aside from the aforementioned static and quasi-static mechanical properties, origami also possesses rich dynamics that can be used to construct prospective engineering devices for impact mitigation and vibration control, which has been studied using experiments, numerical simulations, and modelings~\cite{yasuda2019origami,miyazawa2022topological,huang2022origami,yasuda2020data,zhou2016dynamic}. 

In recent years, the modeling and analysis of dynamical systems via data-driven approaches have grown in popularity because of its equation-free, knowledge-free features~\cite{champion2019data,iten2020discovering,raissi2019physics}. Specifically, origami nonlinear dynamics have been realized with excellent accuracy by machine learning methods based on neural networks~\cite{yasuda2020data}. However, machine learning based on the neural networks is typically computationally intensive and time-consuming at the training stage. Furthermore, the model works as a gigantic black box where the decision processes are difficult to understand and cannot explicitly reveal the underlying physics behind the origami dynamics, despite the fact that a recent study based on the recurrent neural network can mimic the Lyapunov exponent of chaotic origami motions from the hidden layer~\cite{yasuda2020data}.

Dynamic mode decomposition (DMD), as one of the most effective machine learning techniques, has recently been used in many fields for its simplicity and interpretability~\cite{tu2014dynamic}. Given the advantages of DMD and relation between DMD and Koopman mode decomposition, DMD is suitable for characterizing various nonlinear physical and biological systems~\cite{kutz2016dynamic}. Originally, DMD was developed to identify spatiotemporal coherent structures from high-dimensional data in the fluid dynamics community~\cite{schmid2010dynamic}. In addition to fluid dynamics, DMD has been successfully applied to the analysis of biological signals and structural dynamics~\cite{ingabire2021analysis,saito2020data},
as well as more recently to topological metamaterials~\cite{scottarxiv}. In these studies, DMD features fast computation and excellent interpretability by a series of physically coherent structures. However, as for the systems with control, DMD is not capable of uncovering the role that the control plays in the system. In light of this, DMD with control (DMDc), a variant of DMD, is proposed to relate the state and the control of a system, highlighting the importance of the control and improving the accuracy of the DMD model~\cite{proctor2016dynamic,hansen2022swarm}. On the other hand, DMDc often ends up with low accuracy or failure to model several nonlinear systems~\cite{kaiser2018sparse}. Therefore, it would be worthwhile to investigate an improved version of DMD for origami dynamics given the needs of data-driven modeling in this field and inherent deficiencies in the current DMD formulations.

In our work, we introduce a data-driven framework called geometry-informed dynamic mode decomposition (giDMD), which is capable of capturing spatiotemporal dynamics of the Kresling origami structures under excitation. By integrating the geometry information of origami into the DMDc, it becomes possible to not only demonstrate the high accuracy of our approach compared with the DMDc, but also reveal the role of the geometrical parameters in the dynamics, showcasing the interpretability. This approach can also be applied to more complicated Kresling origami structures: origami chain and dual origami structure. In the origami chain, the dynamics across frequencies can be predicted precisely and the frequency of the topological boundary state is identified from the features of giDMD. In the dual origami structure under the excitation in different frequencies, the intrawell periodic motion is predicted accurately. The interwell periodic motion and chaotic dynamics can be identified from the characteristics of giDMD, yet the latter also poses some limitations
(regarding the accuracy of its temporal representation) which are of relevance to consider in future studies. Our approach, which is considerably easier to operate and more interpretable than the machine learning based on the neural networks, offers a general technique to handle the origami dynamics in the presence of geometrical parameters.

\section{Geometry-informed DMD for Kresling origami\label{giDMD}}
\begin{figure}[h]
    \includegraphics[width=0.48\textwidth]{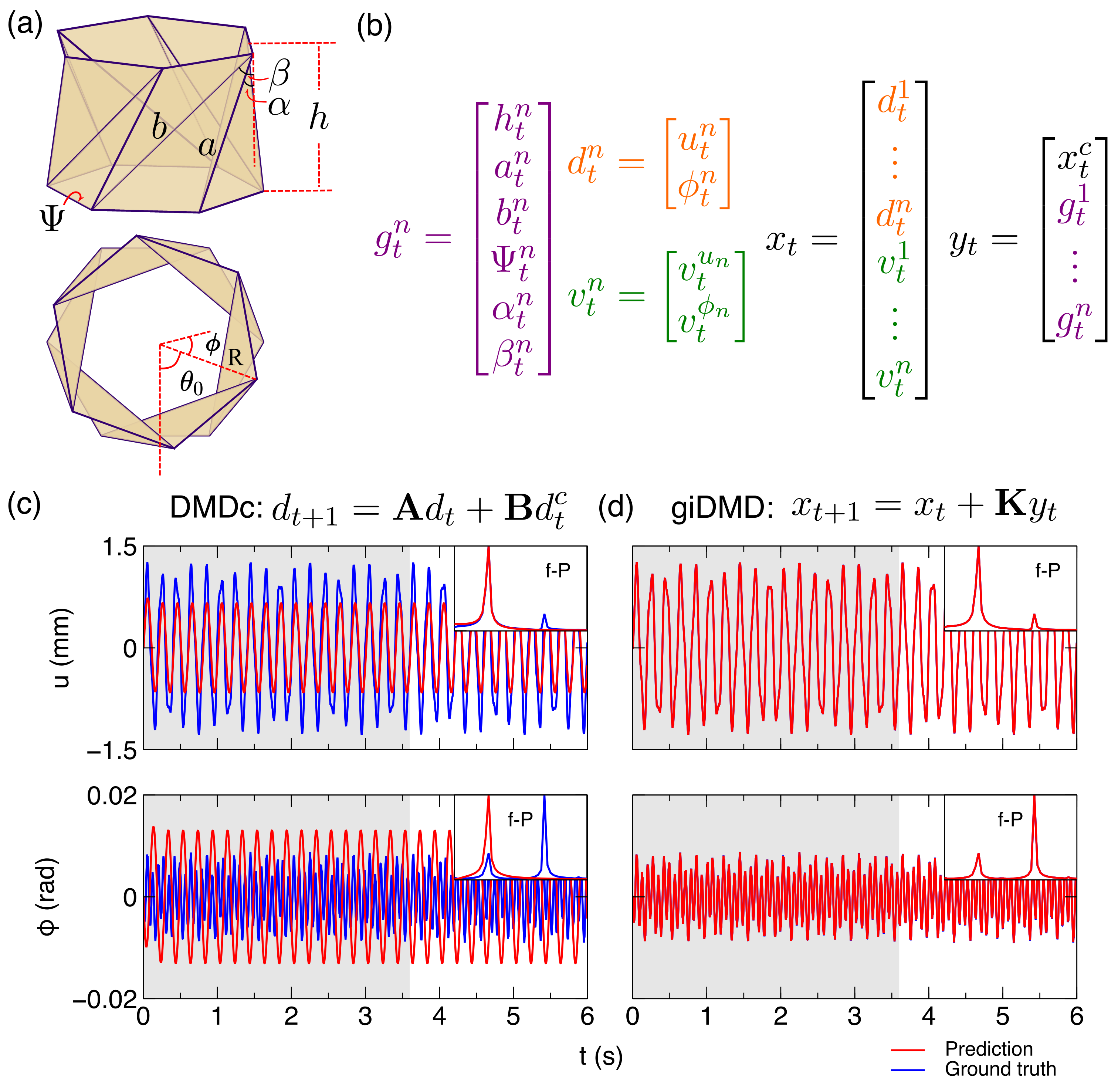}
    \caption{Illustration of giDMD and comparison between DMDc and giDMD.
    (a) The side view and the top view of the single Kresling origami with geometrical parameters.
    (b) The data matrices for the data-driven framework of giDMD. The geometrical parameters $g$ are collected as the height $h$, the length of crease lines $a$, $b$, folding angle $\Psi$, the direction of the crease lines $\alpha$, $\beta$. The state $x$ is taken to include the displacement $p$ and velocity $v$, and the control $y$ is taken to involve the states of the input $x^{c}$ and the geometrical parameters $g$.
    (c) The prediction using DMDc.
    (d) The prediction using giDMD.
    The predictions of axial displacement $u$ and rotational displacement $\phi$ using DMDc are shown in the first and second row, respectively. The gray shaded areas represent the training data. The ground truth and the prediction are shown in blue and red, respectively. The insets in four graphs show the relation between frequency ($f$) and power amplitude ($P$) after the fast Fourier transform.\label{fig:fig1}
    }
\end{figure}

FIG.~\ref{fig:fig1}(a) illustrates the side and top views of the single Kresling origami structure. This single origami structure has two coupled degrees of freedom: translation along the vertical direction $u$ and rotation around vertical direction $\phi$, resulting in the rotation of the top surface following the compression of the origami. This coupled behavior is described by a truss model which is composed of the lumped masses and discs connected by the springs. Appendix~\ref{Appendix A} shows the governing equations of motion for the origami coupled features. The creasing and folding of the origami structure also introduce numerous geometrical variables which are marked in FIG.~\ref{fig:fig1}. $h$, $a$, $b$, $\alpha$, $\beta$ and $\Psi$ represent the height, the length of crease lines, the angles between the crease lines and the vertical direction, and folding angle (the angle between the horizontal plane and facet). During the compression and tension of the origami structure, according to the geometric constraints, these geometrical variables will change as a function of axial and rotational displacement as shown below:
\begin{eqnarray}
    \label{equ:1}
    h=h_{0}-\delta u
\end{eqnarray}
\begin{eqnarray}
    \label{equ:2}
    a=\sqrt{(h_{0}-\delta u)^2+4R^2\sin^2(\frac{\delta \phi}{2}+\frac{\theta_{0}}{2}-\frac{\pi}{2N})}
\end{eqnarray}
\begin{eqnarray}
    \label{equ:3}
    b=\sqrt{(h_{0}-\delta u)^2+4R^2\sin^2(\frac{\delta \phi}{2}+\frac{\theta_{0}}{2}+\frac{\pi}{2N})}
\end{eqnarray}
 \begin{eqnarray}
    \label{equ:4}
    \Psi=\arctan{\frac{h_{0}-\delta u}{R[\cos(\frac{\pi}{N})-\cos(\delta \phi+\theta_{0})]}}
\end{eqnarray}
where the $\delta u$ and $\delta \phi$ represent the differences of axial and rotational displacements between top and bottom surfaces, respectively. $h_{0}$, $\theta_{0}$, $R$ and $N$ are initial height, initial rotation angle representing the chirality, radius and number of vertices of the polygonal cross-section. Accordingly, the angles between the crease lines and the vertical direction can be easily calculated. Under the excitation, axial and rotational displacement will vary over time, hence leading to the variation of geometrical parameters. The vector of geometrical parameters for $n$-th origami element at time $t$ is shown as $g^{n}_{t}$ in FIG.~\ref{fig:fig1}(b).

DMDc is a data-driven approach to model the systems under control, which includes the scenario that the origami structure is excited by an input (the control). DMDc uses a state transition matrix $\bm{A}$ and the control matrix $\bm{B}$ to relate the system displacement $d_{t}$ and the control $d^{c}_{t}$ according to:
\begin{eqnarray}
    \label{equ:5}
    d_{t+1}=\bm{A}d_{t}+\bm{B}d^{c}_{t}
\end{eqnarray}
Depending on whether $\bm{B}$ is known or not, DMDc has different procedures to obtain the model~\cite{proctor2016dynamic}. Herein, we propose a variant that we term geometry-informed DMD based on DMDc. In particular, we refer to geometry-informed DMD~(giDMD) as the DMDc learning framework that integrates underlying knowledge of the origami geometrical variables $g^{n}_{t}$. More concretely, in our present origami setting, 
we create an augmented state of the $n$-th origami $x^{n}_{t}$ composed of the displacement $d^{n}_{t}$ and velocity $v^{n}_{t}$ in both axial and rotational directions, resulting in $x_{t} \in \mathbb{R}^{4n}$. An augmented control $y^{n}_{t}$ is also created by concatenating the augmented state of the control and the vector of geometrical parameters, resulting in $y_{t} \in \mathbb{R}^{l}$, where $l$ is the number of control variables. The vector of geometrical parameters also includes the trigonometric functions (sine and cosine functions) of quantities related to angles, but for simplicity of depicting our approach in a general way, we do not illustrate them in FIG.~\ref{fig:fig1}(b). A similar idea using the augmented state and control variables can also be found in the modeling of swarm dynamics~\cite{hansen2022swarm}. Therefore, the system can be modeled in the form:
\begin{eqnarray}
    \label{equ:6}
    x_{t+1}=x_{t}+\bm{K}y_{t}
\end{eqnarray}
where $\bm{K} \in \mathbb{R}^{4n\times l}$. In giDMD, motivated from the earlier work of~\cite{hansen2022swarm}, the state transition matrix $\bm{A}$ is assumed to be identity ($\bm{A}=\bm{I}$) because we assume that the difference between $x_{t+1}$ and $x_{t}$ ascribes to the contribution of the control. The Equation~(\ref{equ:6}) can be rewritten as:
\begin{eqnarray}
    \label{equ:7}
    x_{t+1}-x_{t}=\bm{K}y_{t}
\end{eqnarray}
Furthermore, Equation~(\ref{equ:7}) can be written in the matrix form:
\begin{eqnarray}
    \label{equ:8}
    \bm{X}'-\bm{X}=\bm{K}\bm{Y}
\end{eqnarray}
\begin{eqnarray}
    \label{equ:9}
    \bm{S}=\bm{K}\bm{Y}
\end{eqnarray}
where $\bm{X}' \in \mathbb{R}^{4n\times (T-1)}$ is one snapshot forward compared with $\bm{X} \in \mathbb{R}^{4n\times (T-1)}$, and $\bm{Y} \in \mathbb{R}^{l\times (T-1)}$ contains $m$ control variables and $(T-1)$ snapshots.

Therefore, the matrix $\bm{K}$ is the key for giDMD modeling, which can be solved by the optimization problem:
\begin{eqnarray}
    \label{equ:10}
    \underset{\bm{K}}{\arg\min} \frac{1}{2}\|\bm{S}-\bm{K}\bm{Y}\|^{2}+\eta R(\bm{K})
\end{eqnarray}
where $R(\cdot)$ is a regularizer that promotes sparsity and $\eta$ is a hyperparameter to determine the strength of the regularization. A sparse matrix $\bm{K}$ can be useful to promote the interpretation of the role geometries play in the origami dynamics and to discover the dominant geometric components. We use the non-convex regularizer, $l_{0}$ norm, in our optimization problem, where sequential thresholded least squares method is performed, which is also widely used in the dynamics discovery~\cite{brunton2016discovering}. $l_{0}$ regularization works by encouraging the coefficients of the model below a certain threshold value to be completely ignored. To be specific, first of all, the standard least squares fitting is performed. Then, a parameter $\eta$ is given to specify the minimum magnitude for coefficients in $\bm{K}$ and all coefficients with magnitude below the threshold are zeroed out. This process of fitting and thresholding is performed until convergence. In this way, a sparse matrix $\bm{K}$ that balances the tradeoff between accuracy and minimizing the number of the control elements can be obtained.

Next, we demonstrate the superiority of our method over the standard DMDc. We simulate the dynamics of a single origami element under the harmonic excitation at the frequency of $5~\mathrm{Hz}$ using the truss model, as detailed in Appendix~\ref{Appendix A}. Due to the coupled behavior of the origami features, the output is the axis displacement $u$ and rotational displacement $\phi$. When the origami dynamics is modeled by DMDc using $60\%$ training data, both predicted displacements $\hat{u}$ and $\hat{\phi}$ exhibit large discrepancies compared with the simulated displacements $u$ and $\phi$ in FIG.~\ref{fig:fig1}(c). The relative errors calculated by $\frac{|u-\hat{u}|_{2}}{|u|_{2}}$ and $\frac{|\phi-\hat{\phi}|_{2}}{|\phi|_{2}}$ are around $50\%$ and $247\%$ for $u$ and $\phi$. In comparison, the significant improvement is achieved by giDMD model, as shown by the agreement between prediction (red lines) and ground truth (blue lines) in FIG.~\ref{fig:fig1}(d). The relative errors for $u$ and $\phi$ are around $0.26\%$ and $1.37\%$. This improvement can also be seen in the spectra (insets in FIG.~\ref{fig:fig1}(c) and FIG.~\ref{fig:fig1}(d)) given by the fast Fourier transform. Under the excitation of $5~\mathrm{Hz}$, there is another peak other than the excitation frequency in the spectrum because of the nonlinearity. In stark contrast with the DMDc, giDMD can predict two peaks in the spectrum, implying that giDMD can fully capture the origami dynamics. The 
drastic improvement in accuracy in both time and frequency domain demonstrates the effectiveness of giDMD towards modeling origami dynamics.

\section{Topological boundary states in the origami chain\label{Topo}}
\begin{figure}[h!]
    \includegraphics[width=0.5\textwidth]{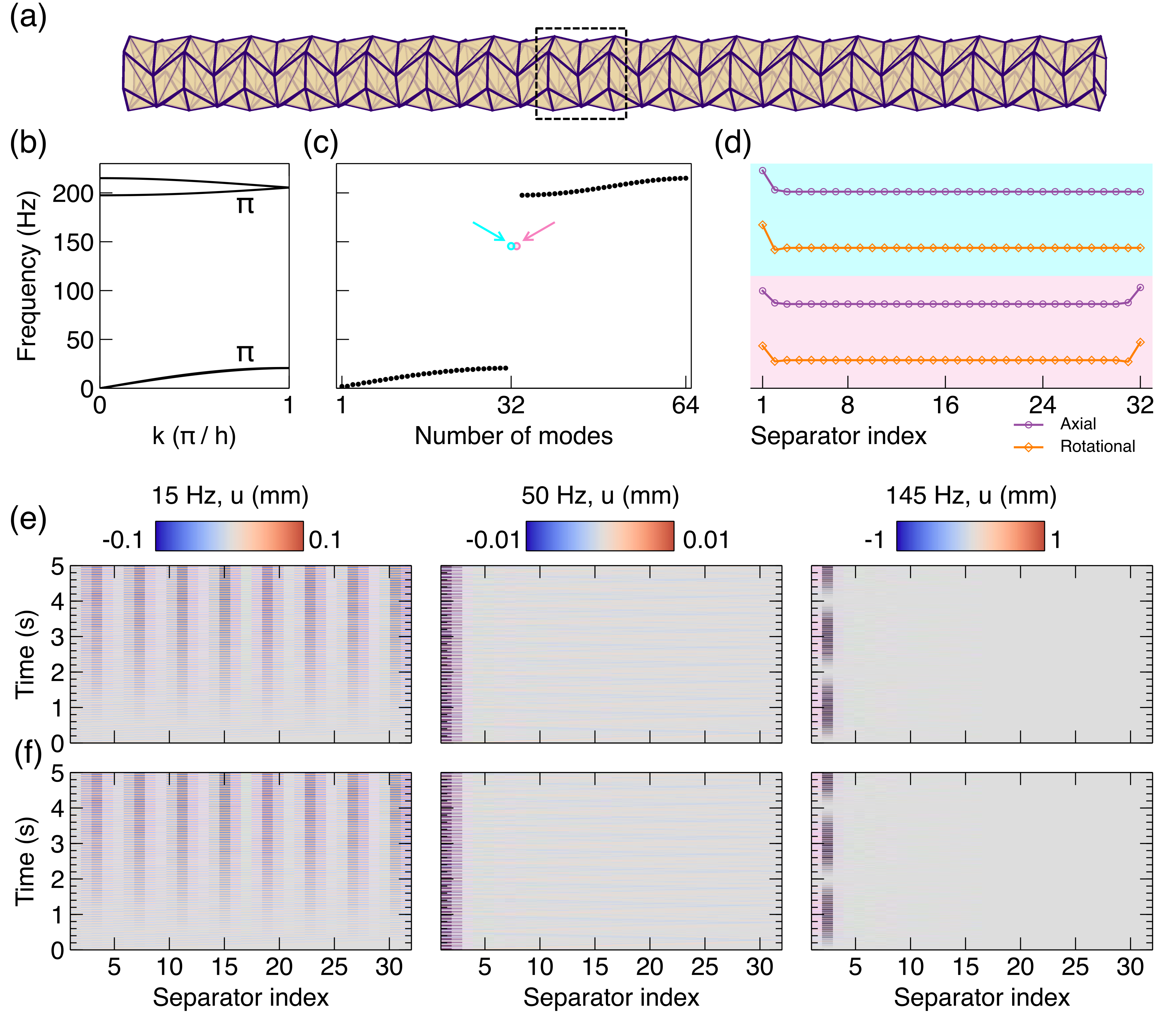}
    \caption{Application of giDMD on the origami chain.
    (a) The schematic of the origami chain composed of Kresling origami with alternating chirality. The design parameters of origami are $h_{0}=30~\mathrm{mm}$, $\theta_{0}=\pm70^{\circ}$, $R=36~\mathrm{mm}$ where $h_{0}$, $\theta_{0}$, and $R$ are the initial height, initial rotational angle, and radius of the cross-section, respectively.
    (b) The band structure calculated by the unit cell enclosed by black dashed line in (a). The Zak phases of the lower band and upper band are $\pi$.
    (c) The modes of the supercell formed by $16$ unit cells. The black dots represent the bulk band. The cyan and pink dots represent the topological boundary states within the band gap.
    (d) The axial and rotational modes of the topological boundary states.
    (e) The simulated axial displacement along the origami chain at different frequencies in the pass band~($15~\mathrm{Hz}$), stop band~($50~\mathrm{Hz}$) and topological boundary states~($145~\mathrm{Hz}$).
    (f) The corresponding predicted axial displacement using $60$\% training data~(from $0~\mathrm{s}$ to $3~\mathrm{s}$).
    \label{fig:fig2}}
\end{figure}

The effectiveness of our method has emerged in the single origami structure above. Next, we further demonstrate this method on more complex origami structures. The first example we demonstrate our method on is the elastic topological metamaterials built by the origami structure. The elastic topological metamaterials are inspired by the electronic topological insulators, where the vibration is isolated in the bulk yet propagating along the boundary or surface. The topologically protected defect-immune wave propagation in the elastic topological metamaterials has been attracting significant attention~\cite{ma2019topological,li2019valley,liu2018tunable,li2021topological,vila2017observation,li2018observation}. Recently, the origami structures have been used as building blocks to construct topological metamaterials~\cite{miyazawa2022topological}. We adopt the existing design to construct the origami chain shown in FIG.~\ref{fig:fig2}(a). The unit cell enclosed by black dashed line contains two origamis with opposite chirality represented by initial rotational angle $\theta_{0}$.
\par
The band structure can be calculated by the linearized truss model (Appendix~\ref{Appendix A}) based on the small amplitude approximation after the application of periodic boundary condition (Bloch's theorem). As shown in FIG.~\ref{fig:fig2}(b), four bands appear in the first Brillouin zone. Two lower bands nearly overlap, and two upper bands are crossing to form a degeneracy point at the edge of the first Brillouin zone. To characterize the topology of this system, we calculate the topological invariant Zak phase for the 1D system. Because the bands are degenerate at the edge of the first Brillouin zone, the topological invariant Zak phase $\varphi$ for $q$-th and $(q+1)$-th bands is calculated through the Wilson-loop eigenvalues as indicated below:
\begin{eqnarray}
    \label{equ:11}
    &&\varphi_{q,q+1}=\\ \nonumber
    &&-\sum_{k=0}^{k-1} \operatorname{Im} \ln{
    \begin{vmatrix}
        \langle U^{q}_{k} | U^{q}_{k} \rangle & \langle U^{q}_{k} | U^{q+1}_{k+1} \rangle \\
        \langle U^{q+1}_{k} | U^{q}_{k+1} \rangle & \langle U^{q+1}_{k} | U^{q+1}_{k+1} \rangle
    \end{vmatrix}
    }
\end{eqnarray}
where $U^{q}_{k}$ denotes the mode at the Bloch wave vector $k$ for the $q$-th band~\cite{wang2019band}. This produces $\varphi_{1,2}=\varphi_{3,4}=\pi$ which are marked in FIG.~\ref{fig:fig2}(b). The nonzero topological invariant of the lower bands also ensures the topologically nontrivial band gap between the lower bands and the upper bands.

According to the bulk-edge correspondence, the topological boundary state will emerge within the band gap in the truncated origami chain. In FIG.~\ref{fig:fig2}(c), the eigenmodes of the truncated origami chain with $16$ unit cells is calculated. As expected, two degenerate modes appear within the band gap. After checking the axial and rotational mode shape at $k_{x}=0$ in FIG.~\ref{fig:fig2}(d), we confirm that these two degenerate modes are topological boundary states due to the concentrated displacement at the boundary.

Then, the wave propagation along the origami chain is calculated by the truss model~(Appendix~\ref{Appendix A}) using the small-amplitude approximation~(i.e., in the linear regime). The axial displacements of each separator over time are illustrated in FIG.~\ref{fig:fig2}(e) under different excitation frequencies, represented by the frequencies in pass band~($15~\mathrm{Hz}$), stop band~($50~\mathrm{Hz}$) and at the topological boundary state~($145~\mathrm{Hz}$). At the frequency of $15~\mathrm{Hz}$, the vibration can affect the whole origami chain, while the vibration is localized at the boundary at the frequency of $50~\mathrm{Hz}$, as expected for such a band gap frequency. At the frequency of $145~\mathrm{Hz}$, where the topological boundary state emerges, the vibration is also localized at the boundary but with larger amplitude, representing the excitation of the corresponding eigenmode.

The giDMD model for the origami chain is put forward using $60\%$ training data~(from $0~\mathrm{s}$ to $3~\mathrm{s}$). The prediction and ground truth agree excellently, as evidenced by the same patterns in FIG.~\ref{fig:fig2}(e) and FIG.~\ref{fig:fig2}(f), and relative error smaller than $0.1\%$ in general. The highly accurate modeling of origami chain across frequencies confirms the ability of giDMD to model complex structures in the nearly linear regime of the corresponding dynamics.

\begin{figure}[h]
    \includegraphics[width=0.5\textwidth]{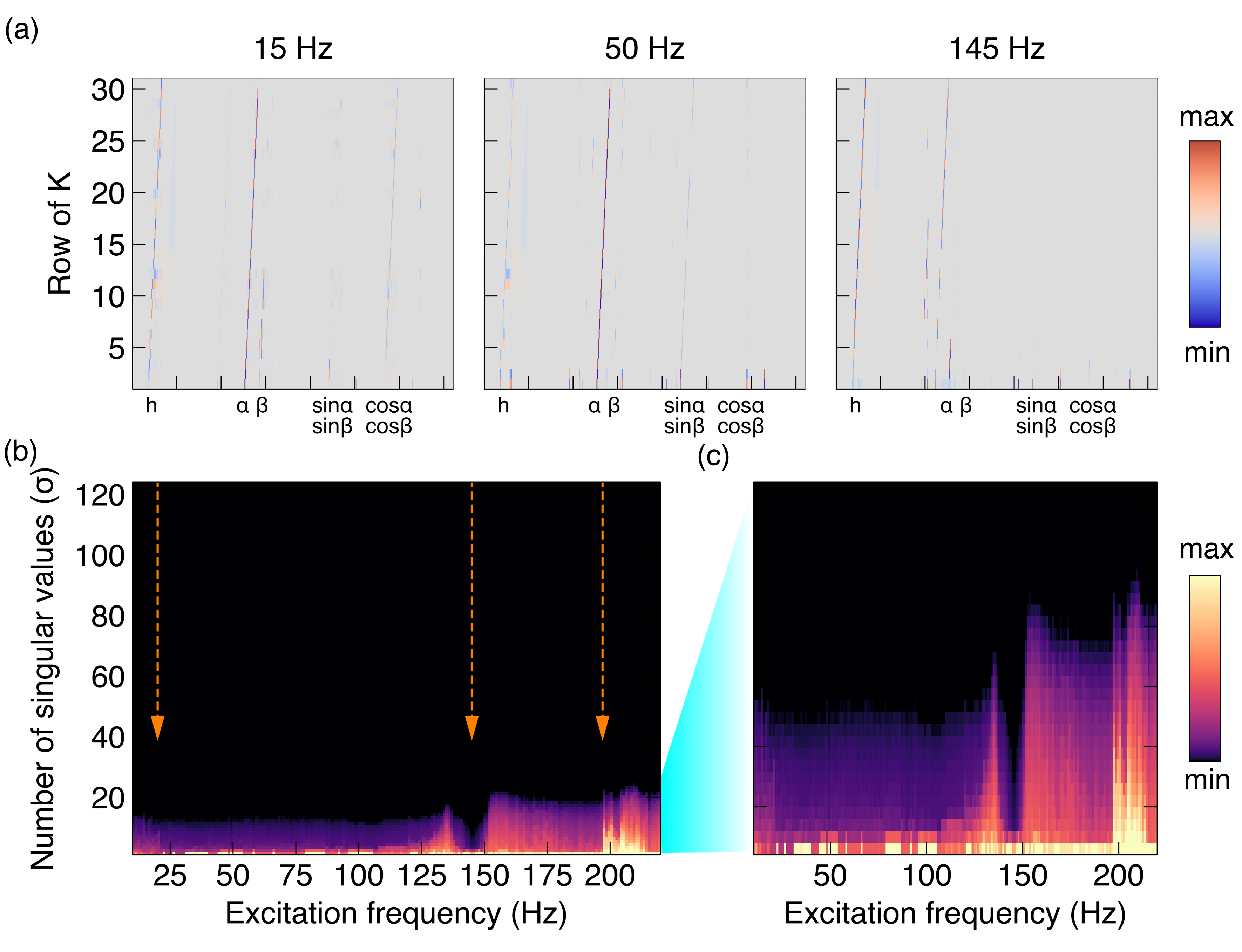}
    \caption{Identification of topological boundary states in origami chain.
    (a) The rows of $\bm{K}$ matrix corresponding to the axial displacements of the origami chain at $15~\mathrm{Hz}$, $50~\mathrm{Hz}$ and $145~\mathrm{Hz}$. The corresponding control variables are shown below.
    (b) The singular value spectrum of $\bm{K}$ matrix under different excitation frequencies from $10$ to $220~\mathrm{Hz}$. The orange dash arrows on two sides indicate the range of the bandgap, and one in the middle indicates the frequency of topological boundary states.
    (c) The zoom-in view by using $32$ singular values.
    \label{fig:fig3}}
\end{figure}

We take the first $31$ rows of the $\bm{K}$ matrix responsible for the calculation of the axial displacement of each separator in the origami chain. As shown in FIG.~\ref{fig:fig3}(a), the $\bm{K}$ matrix is fairly sparse for all three cases, where nonzero values are concentrated near the position indicating the information of height $h$, angles between crease lines and vertical direction $\alpha$, $\beta$, and corresponding sine functions $\sin \alpha$, $\sin \beta$. It suggests that these geometrical variables are mainly important and responsible for the axial displacement of origami chain. Especially, at the frequency of the topological boundary state~($145~\mathrm{Hz}$), the values for $\sin \alpha$, $\sin \beta$ are smaller than those of the other two frequencies, suggesting less importance of these two geometrical variables in forming topological boundary states. This may result from the strong localization of vibrations near the boundary, causing most geometrical variables to remain almost unchanged. Furthermore, the nonzero values shift linearly with the increase of the row of $\bm{K}$, indicating that the displacement of each separator is related to geometrical variables of nearest origami elements.

Similar to the DMD where the eigenvalue of the linear operator can characterize the system dynamics, we use the singular 
values of the $\bm{K}$ matrix to identify the `state' of the system. The singular value spectra of the $\bm{K}$ matrix from $10~\mathrm{Hz}$ to $220~\mathrm{Hz}$ are shown in FIG.~\ref{fig:fig3}(b) to illustrate the difference. In the region of the pass band~($f~\le~20~\mathrm{Hz}$, $f~\ge~198~\mathrm{Hz}$), the first several singular values are 
significantly larger than those in the region of stop band~($20~\mathrm{Hz}\le~f~\le~198~\mathrm{Hz}$). At the frequency of topological boundary states, the singular values are particularly smaller than any of other frequencies as detailed in zoom-in view in FIG.~\ref{fig:fig3}(c). The smaller singular values at the frequency of topological boundary state may be induced by the absence of contribution from $\sin \alpha$, $\sin \beta$ compared with other frequencies. We also notice that there are smaller singular values around $145~\mathrm{Hz}$ but the smallest one is at $145~\mathrm{Hz}$. The singular values are large outside $145~\mathrm{Hz}$ and the ones within the band gap region~($20~\mathrm{Hz}\le~f~\le~198~\mathrm{Hz}$) are used to describe the localized states near the boundary due to the band gap.

\section{Dynamic motion in the dual origami structure\label{dual}}
\begin{figure}[h]
    \centering
    \includegraphics[width=0.5\textwidth]{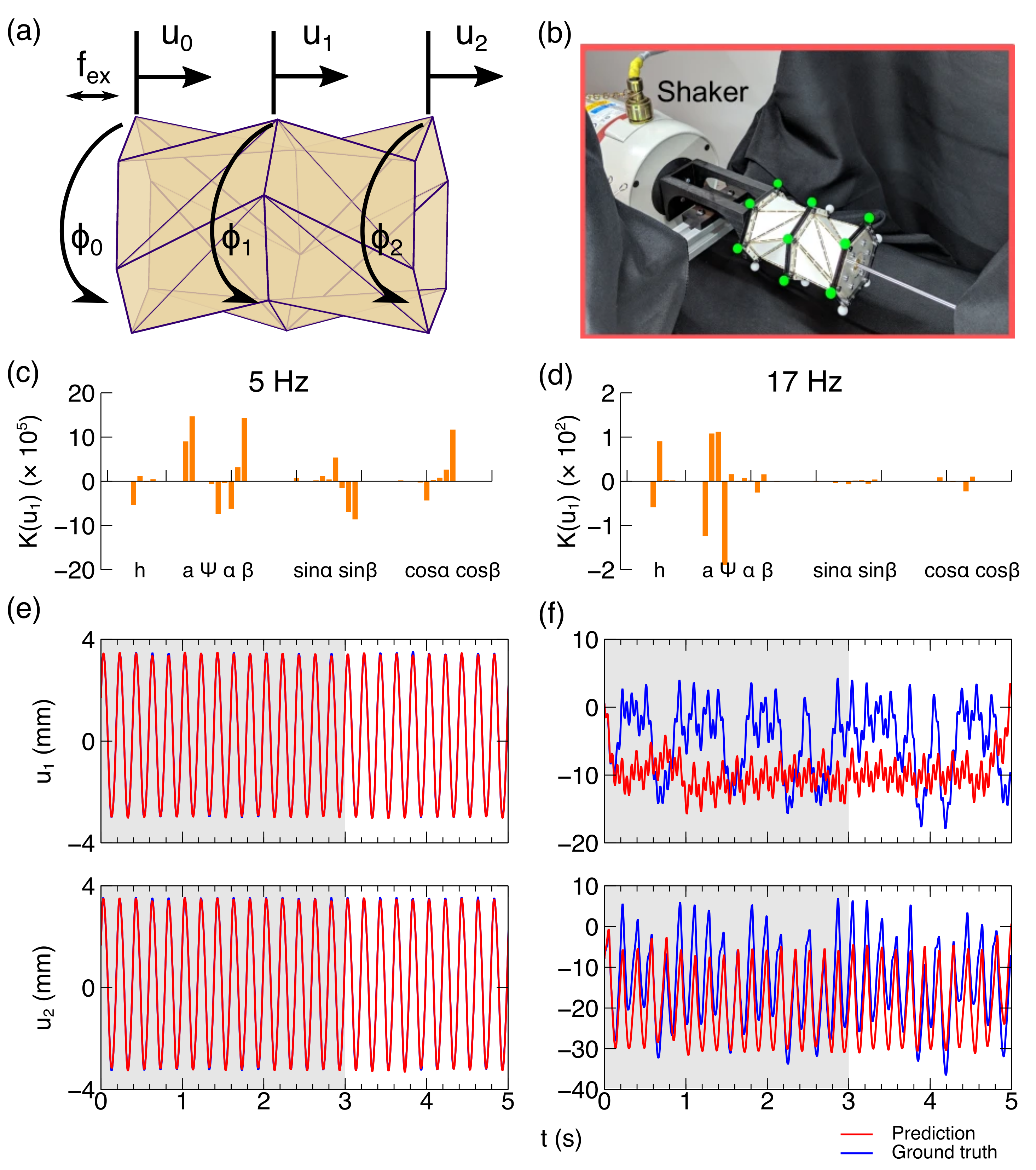}
    \caption{Application of giDMD on the dual origami structure.
    (a) The schematic of the dual origami structure with opposite chirality.
    (b) The image of experimental setup for the vibration test from the reference~\cite{yasuda2020data}. Copyright 2020 Springer Nature. The design parameters of origami are $h_{0}=50~\mathrm{mm}$, $\theta_{0}=\pm70^{\circ}$, $R=36~\mathrm{mm}$ where $h_{0}$, $\theta_{0}$, and $R$ are the initial height, initial rotational angle, and radius of the cross-section, respectively. The input excitation is applied by a shaker, and the motions are captured by two action cameras and quantified by the digital image correlation program.
    (c)(d) The row of $\bm{K}$ matrix corresponding to the axial displacement $u_{1}$ at the frequency of $5~\mathrm{Hz}$ and $17~\mathrm{Hz}$, respectively. The corresponding control variables are shown below.
    (e) The prediction of axial displacement $u_{1}$ and $u_{2}$ using giDMD at the frequency of $5~\mathrm{Hz}$.
    (f) The prediction of axial displacement $u_{1}$ and $u_{2}$ using giDMD at the frequency of $17~\mathrm{Hz}$.
    The gray shaded areas represent the training data. The ground truth is shown in blue and the prediction is shown in red.
    \label{fig:fig4}}
\end{figure}

The next example is the dual origami structure exhibiting rich dynamic motions under excitation. This structure is composed of two bistable origami elements with opposite chirality as shown in FIG.~\ref{fig:fig4}(a), in which $u_{0}$ and $\phi_{0}$ are the input axial and rotational excitations at the first separator, and $u_{1}$, $\phi_{1}$ and $u_{2}$, $\phi_{2}$ are the axial and rotational displacements at the second and third separators, respectively. In the previously 
reported experimental study~\cite{yasuda2020data}, under the excitation in different frequencies, the structure will 
feature periodic motion~($5-9~\mathrm{Hz}$, $14-16~\mathrm{Hz}$, $18~\mathrm{Hz}$, $23-24~\mathrm{Hz}$) and chaotic motion~($10-13~\mathrm{Hz}$, $17~\mathrm{Hz}$, $19-22~\mathrm{Hz}$). Note that due to combination of two bistable origamis, the periodic motion can be further identified as intrawell periodic motion~($5-9~\mathrm{Hz}$, $18~\mathrm{Hz}$, $23-24~\mathrm{Hz}$), and interwell periodic motion~($14-16~\mathrm{Hz}$). The intra- and inter- well behaviors are further elaborated in the reference~\cite{yasuda2020data}. The experimental setup is shown in FIG.~\ref{fig:fig4}(b) and detailed in Appendix~\ref{Appendix A}. Although this structure is simpler than the previous example, it is under the large-amplitude excitation region and hence displays the corresponding hallmarks of nonlinear dynamics, including the above mentioned chaotic motion.

To demonstrate our method, we choose two typical cases where periodic motion and chaotic motion are represented by $5~\mathrm{Hz}$ and $17~\mathrm{Hz}$, respectively. After obtaining the $\bm{K}$ matrix from the sparse regression, we check the first row of the $\bm{K}$ matrix responsible for the calculation of axis displacement $u_{1}$. The case for $5~\mathrm{Hz}$ is shown in FIG.~\ref{fig:fig4}(c). The first row of the $\bm{K}$ matrix is sparse so that nonzero values only appear in several positions standing for the information of $h$, $a$, $\Psi$, $\alpha$, $\beta$, sine and cosine functions of $\alpha$ and $\beta$, suggesting that these geometrical variables are responsible for the axial motion of the second separator. In comparison, at the frequency of $17~\mathrm{Hz}$, nonzero values emerge in similar positions but with much smaller values in FIG.~\ref{fig:fig4}(d). Besides, the sine and cosine functions of $\alpha$ and $\beta$ are contributing less to the axial motion of the second separator due to the smaller values than other geometrical variables. Compared with the origami chain where the length of crease lines and the folding angle do not significantly contribute to the axial displacement, in the present case, the length of crease lines and folding angle are major contributions to the axial displacement for the dual origami structure. The possible reason can be the smaller amplitude excitation~(i.e., linear regime) of the origami chain, leading to the representation of $a$, $b$ and $\Psi$ by the function of $h$. After sparse regression, the contributions of $a$, $b$ and $\Psi$ can be approximated by $h$, whereas large-amplitude excitation~(nonlinear regime) in dual origami structure results in the failure to represent $a$, $b$ and $\Psi$ using $h$.

We then predict the axial displacements $u_{1}$ and $u_{2}$ at the frequency of $5~\mathrm{Hz}$ and $17~\mathrm{Hz}$. As shown in FIG.~\ref{fig:fig4}(e), the prediction and ground truth agree excellently with $60\%$ training data, as evidenced by the overlap of red curves and blue curves. In stark contrast, the predictions at the frequency of $17~\mathrm{Hz}$ do not match with the ground truth in a precise manner~(FIG.~\ref{fig:fig4}(f)). Although the sudden change in the chaotic motion cannot be predicted, the predicted results still show clear differences from the periodic motion~(FIG.~\ref{fig:fig4}(e)), i.e., they carry the blueprint of the observed chaotic dynamics. Therein, $u_{1}$ is more chaotic than $u_{2}$ and thus the prediction is better in $u_{2}$. Note that the prediction results of intrawell periodic motion at other frequencies are similar to the FIG.~\ref{fig:fig4}(e), while the ones of chaotic motion at other frequencies are similar to the FIG.~\ref{fig:fig4}(f).

\begin{figure}[h]
    \centering
    \includegraphics[width=0.5\textwidth]{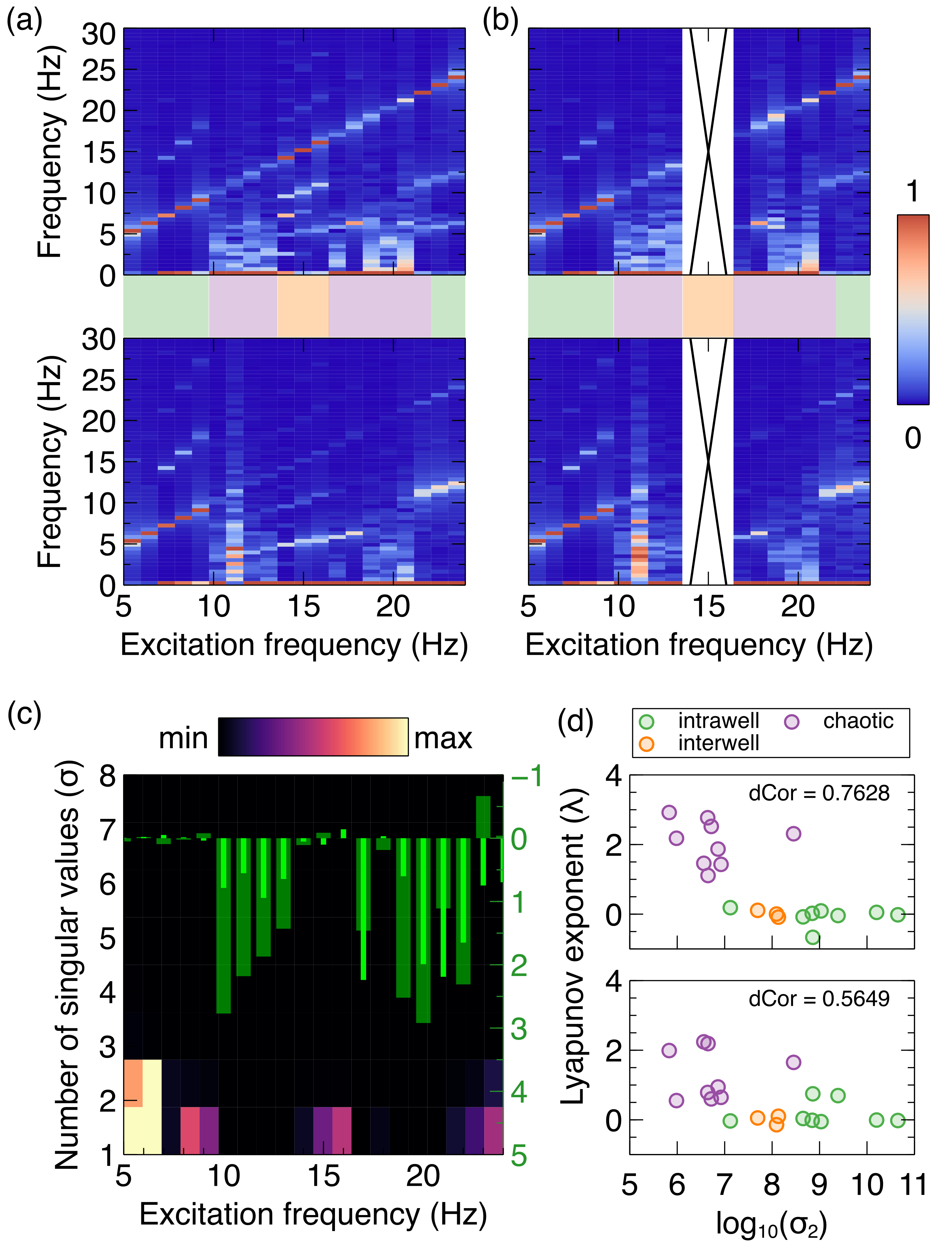}
    \caption{Identification of intrawell periodic, interwell periodic and chaotic motion.
    (a) The spectrum analysis of the experimental axial displacements $u_{1}$ (the first panel) and $u_{2}$ (the second panel) under different excitation frequencies from $5$ to $24~\mathrm{Hz}$.
    (b) The corresponding spectral analysis based on the predicted axial displacements $u_{1}$ (the first panel) and $u_{2}$ (the second panel). The white area with a cross represents the failure of the prediction.
    The green, orange and purple shaded areas represent the frequencies of intrawell periodic motion, interwell periodic motion and chaotic motion.
    (c) The singular value spectrum of $\bm{K}$ matrix under different excitation frequencies from $5$ to $24~\mathrm{Hz}$. The Lyapunov exponents of $u_{1}$ and $u_{2}$ under different excitation frequencies are also shown in dark green and bright green, respectively.
    (d) The relation between the second singular values $\sigma_{2}$ and the Lyapunov exponents $\lambda$ under different excitation frequencies. The intrawell periodic motion, interwell periodic motion and chaotic motion are represented in green, orange and purple circles. The distance correlation between logarithm-scale second singular values and Lyapunov exponents are shown in both graphs.
    \label{fig:fig5}}
\end{figure}

The spectral analysis is conducted to show the frequency response under the excitation in different frequencies. In FIG.~\ref{fig:fig5}(a), the spectra for $u_{1}$ and $u_{2}$ in different frequencies are obtained by the fast Fourier transform of the experimental data after normalization. It is obvious that along the diagonal direction of each panel, there are responses at the same frequency as the excitation frequency. Apart from that, the periodic motion and the chaotic motion can be clearly identified from the increase of lower frequency components.
The `state' of the structure is also marked in FIG.~\ref{fig:fig5}(a) and FIG.~\ref{fig:fig5}(b), where green, orange and purple shaded areas indicate the intrawell periodic motion, interwell periodic motion and chaotic motion. FIG.~\ref{fig:fig5}(b) illustrates the spectra for $u_{1}$ and $u_{2}$ based on the prediction results. It can be seen that the spectra in the region of intrawell periodic motion~($5-9~\mathrm{Hz}$, $18~\mathrm{Hz}$, $23-24~\mathrm{Hz}$) agree with the ground truth excellently. However, giDMD fails to predict the interwell periodic motion~($14-16~\mathrm{Hz}$), where the predictions blow up to infinity at some time point, resulting in the failure of spectral analysis. The interwell periodic motion describes the scenario where $u_{1}$ and $u_{2}$ feature a large difference, as well as $\phi_{1}$ and $\phi_{2}$, i.e., in the
regime of significant axial and rotational differences across origami separators. In such case, the blow-up phenomena in the prediction are due to the unstable dynamical system produced by the giDMD model. Understanding how to extend the giDMD approach to capture such a scenario remains an interesting open question for further investigation. In the region of chaotic motion~($10-13~\mathrm{Hz}$, $17~\mathrm{Hz}$, $19-22~\mathrm{Hz}$), the spectra from the predicted results show qualitative agreement with the ground truth featured by the increased magnitude in the low frequency range, despite the failure to capture the intermittent variation in chaotic motion as shown in FIG.~\ref{fig:fig5}(f). These findings imply that giDMD can predict the chaotic motion to some extent.

Similar to the section~\ref{Topo}, we study the $\bm{K}$ matrix to identify the `state' of the dual origami structure. FIG.~\ref{fig:fig5}(c) shows the singular values $\sigma$ of the $\bm{K}$ matrix at different frequencies, where different singular value spectra appear in different frequencies. We notice that in the frequency range corresponding to the chaotic motion, the singular values are particularly smaller than others. To characterize the origami dynamics, the Lyapunov exponents for $u_{1}$ and $u_{2}$ are shown by bright green and dark green in FIG.~\ref{fig:fig5}(c), calculated by the Rosenstein’s method~\cite{rosenstein1993practical}. 
As is well-known, the larger Lyapunov exponent suggests more chaotic dynamics. From our results, we can deduce that smaller singular values correspond to larger Lyapunov exponent. In FIG.~\ref{fig:fig5}(d), we choose the second singular value of each $\bm{K}$ matrix $\sigma_{2}$ to compare with Lyapunov exponents for $u_{1}$~(top panel) and $u_{2}$ (bottom panel). It is noticed that there are clear boundaries among chaotic motion~(smallest values), interwell periodic motion and intrawell periodic motion~(largest values). Besides, there is the fairly strong relation characterized by the distance correlation~(quotient of the distance covariance and the product of the distance standard deviations) between logarithm-scaled second singular values and Lyapunov exponents, resulting in $dCor=0.7628$ and $dCor=0.5649$ for $u_{1}$ and $u_{2}$. The correlation between the singular values and Lyapunov exponent of the dynamics system has emerged in several studies and has been used to calculate Lyapunov exponent~\cite{geist1990comparison,barreira2017lyapunov}. Therefore, from the singular values of $\bm{K}$, one can identify the nature of the dynamics of the dual origami structure.

\section{Discussion and Conclusion\label{conclusion}}
The $\bm{K}$ matrix serves as the control matrix in the giDMD, which will convert the control $\bm{Y}$ to $\bm{S}$, expressed by $\bm{S}=\bm{X}-\bm{X}'$. Since the state $\bm{X}$ corresponds to the displacement and velocity of each separator at each time step, $\bm{S}$ contains the information of velocity and acceleration. Essentially, $\bm{K}$ is found to describe the velocity and acceleration using geometrical parameters in origami structures. The sparsity of $\bm{K}$ matrix shows that the velocity and acceleration of the origami structure are only related to a few geometrical parameters. This can also be implied by the governing equation of motion of origami as shown in the Appendix~\ref{Appendix A} and derived equations. These governing equations of motion of origami contain the combination of geometrical parameters explicitly or implicitly. This is the reason why giDMD works well in origami dynamics and is only related to the several geometrical variables, resulting in the sparse $\bm{K}$ matrix. 
In a sense, one can argue that the material in Appendix~\ref{Appendix A} represents the traditional modeling approach (based on constitutive laws) towards the origami system of interest, while the giDMD represents a modern, data-driven variant thereof. Indeed, our method, from the idea of dynamics discovery using sparse regression, is quite similar to the sparse identification of nonlinear dynamics~(SINDy)~\cite{brunton2016discovering}. The sparse $\bm{K}$ in giDMD and the sparse $\bm{\Xi}$ in SINDy both serve to select the active terms in the library (control variables in giDMD and candidate functions in SINDy). However, there are some differences. In SINDy, there exists an extensive library of candidate functions that is provided as a possible basis to represent the dynamics, while $\bm{Y}$ in giDMD only contains the origami geometrical variables. Moreover, the exact velocity and acceleration of the system are used in the SINDy, while $\bm{S}$ in giDMD, motivated by the DMDc approach, is composed of the difference of displacement and the difference of velocity, which is different from the implementation of SINDy. For the consistency of the narrative, we fix the training ratio to be $60\%$ in the main text, but we also show the giDMD modeling in other training ratios in Appendix~\ref{Appendix B} and Appendix~\ref{Appendix C} including error analysis and identification of states of origami structures.

In this work, we propose the method of giDMD to learn the origami dynamics that gives rise to origami coupled motions from pure observation data (experiment data and simulation data). We show the better performance of giDMD characterized by the prediction accuracy compared with DMDc. In the two example origami structures we apply our method on, giDMD can not only predict the origami behaviors under different frequencies, but also highlight an ability to identify the `state' of origami structures. Furthermore, giDMD offers the insights into the importance of geometries in the governing motion law of origami. giDMD has better performance in the linear (or nearly linear) region of origami dynamics, while there are still some challenges in describing chaotic dynamics within the nonlinear regime. In the latter, giDMD starts to fail in connection to the goal of accurate prediction, but it can still help with the identification of the origami state. The giDMD provides a substantial capability to model origami dynamics in an efficient~(fast computation) and interpretable way. Although we primarily demonstrate our method on two specific structures, the studied origami structures can be more complex across scales and not limited to the Kresling origami in the future work. Indeed, the further extension of giDMD to other such geometrically nonlinear systems is a topic currently under active investigation and relevant results will be reported in due course.

\begin{acknowledgments}
We are grateful for the support from the U.S. National Science Foundation: S.L., Y.M., K.Y. and J.Y. (Grant No. 1933729), P.G.K. (Grant No. DMS 2204702). J.Y. acknowledges the support from the SNU-IAMD and the Brain Pool Plus program funded by the Ministry of Science and ICT through the National Research Foundation of Korea (0420-20220160). K.Y. is supported by Funai Foundation for Information Technology.
\end{acknowledgments}

\appendix

\section{Data acquisition}
\label{Appendix A}
\subsection{Experimental data}
In the example of the dual origami structure, the raw 
experimental data was used for the giDMD model. The 
experimental data and experimental procedures were detailed in the reference~\cite{yasuda2020data}. The dynamic test was conducted on a dual triangulated cylindrical origami structure. The design parameters of origami were $h_{0}=50~\mathrm{mm}$, $\theta_{0}=\pm70^{\circ}$, $R=36~\mathrm{mm}$ which are the initial height, initial rotational angle, and radius of the cross-section, respectively. The left origami is in negative chirality ($\theta_{0}=-70^{\circ}$) and the right origami is in positive chirality ($\theta_{0}=+70^{\circ}$). The left end of dual origami structure was connected to the shaker (LDS V406 M4-CE, Br\"{u}el \& Kjær) which applied the single-frequency harmonic excitation to the system. The motion of spherical markers attached to the separators was captured by two action cameras in $240$ frames per second and quantified by digital image correlation. The triangulation method was used to find the three-dimensional coordinates of the spherical markers, resulting in the axial displacement $u$ and rotational displacement $\phi$ of each separator.
\subsection{Simulation data}
In the example of single origami and origami chain, the simulation data was used for the giDMD model. The design parameters and mechanical parameters used in the numerical simulation are shown in TABLE~\ref{table:1} and TABLE~\ref{table:2}, respectively.
\begin{table}[h]
    \centering
    \caption{Design parameters for single origami and origami chain}
    \label{table:1}
    \begin{tabular}{c|c|c|c|c|c}
        $m~\mathrm{(kg)}$ & $j~\mathrm{(kg\cdot m^{2})}$ & $N$ & $h_{0}~\mathrm{(m)}$ & $R~\mathrm{(m)}$ & $\theta_{0}~\mathrm{(^{\circ})}$\\
        \hline
        $58.8\times10^{-3}$ & $6.77\times10^{-5}$ & $6$ & $30\times10^{-3}$ & $36\times10^{-3}$ & $\pm70$
    \end{tabular}
\end{table}
\begin{table}[h]
    \centering
    \caption{Mechanical parameters for single origami and origami chain}
    \label{table:2}
    \begin{tabular}{c|c|c}
         $k_{a}~\mathrm{(N\cdot m^{-1})}$ & $k_{b}~\mathrm{(N\cdot m^{-1})}$ & $k_{\Psi}~\mathrm{(N\cdot m\cdot rad^{-1})}$\\
        \hline
         $6055$ & $3743$ & $7.277\times10^{-3}$\\
    \end{tabular}
\end{table}

The mechanical parameters $k_{a}$, $k_{b}$ represent the axial spring constant along the crease $a$, $b$, and $k_{\Psi}$ represents the torsional spring constant along the bottom crease.

The simulations of the single origami and origami chain under excitation in different frequencies are conducted based on the truss model regarding Kresling unit cell as inter-polygonal spring connecting separators with mass $m$ and rotational inertia $j$~\cite{yasuda2017origami,yasuda2019origami,miyazawa2022topological}. The equation of motion of the separator can be expressed as:
\begin{subequations}
\begin{eqnarray}
    \label{equ:12}
    m_{n}\Ddot{u}_{n}+F_{n}(\delta u_{n},\delta \phi_{n})-F_{n-1}(\delta u_{n-1},\delta \phi_{n-1})=0
\end{eqnarray}
\begin{eqnarray}
    \label{equ:13}
    j_{n}\Ddot{\phi}_{n}+T_{n}(\delta u_{n},\delta \phi_{n})-T_{n-1}(\delta u_{n-1},\delta \phi_{n-1})=0
\end{eqnarray}
\end{subequations}
where $\delta u_{n}=u_{n}-u_{n+1}$ and $\delta \phi_{n}=\phi_{n}-\phi_{n+1}$. The subscript $n$ denotes the $n$-th separator. The force and torque in Equation~(\ref{equ:12}) and Equation~(\ref{equ:13}) can be further expanded by the summation of the contribution from each spring as below:
\begin{widetext}
\begin{subequations}
\begin{eqnarray}
    \label{equ:14}
    F_{n}(\delta u_{n},\delta \phi_{n})=&&k_{an}N(\delta u_{n}-h_{0})(1-\frac{a_{0}}{a_{n}})
    +k_{bn}N(\delta u_{n}-h_{0})(1-\frac{b_{0}}{b_{n}})\nonumber\\
    &&+2k_{\Psi n}NRh_{0}(\Psi_{0}-\Psi_{n})
    \frac{\cos \frac{\pi}{N}-\cos(\delta \phi_{n}+\theta_{0})}{R^{2}[\cos \frac{\pi}{N}-\cos(\delta \phi_{n}+\theta_{0})]^{2}+(h_{0}-\delta u_{n})^{2}}
\end{eqnarray}
\begin{eqnarray}
    \label{equ:15}
    T_{n}(\delta u_{n},\delta \phi_{n})=&&k_{an}NR^{2}\sin(\delta \phi_{n}+\theta_{0}-\frac{\pi}{N})(1-\frac{a_{0}}{a_{n}})
    +k_{bn}NR^{2}\sin(\delta \phi_{n}+\theta_{0}+\frac{\pi}{N})(1-\frac{b_{0}}{b_{n}})\nonumber\\
    &&+2k_{\Psi n}NRh_{0}(\Psi_{0}-\Psi_{n})
    \frac{\sin(\delta \phi_{n}+\theta_{0})(h_{0}-\delta u_{n})}{R^2[\cos \frac{\pi}{N}-\cos(\delta \phi_{n}+\theta_{0})]^{2}+(h_{0}-\delta u_{n})^{2}}
\end{eqnarray}
\end{subequations}
\end{widetext}
where the parameters with subscript $0$ denote the initial origami configuration.
\par
The calculation of the band structure is based on the linearized truss model of the unit cell in order to represent the small amplitude excitation. The unit cell contains two origami elements whose front separators denote the sites ($1$) and ($2$). The linearized equations of motion can be expressed as:
\begin{widetext}
\begin{subequations}
\begin{eqnarray}
    \label{equ:16}
    m_{n}\Ddot{u}^{(1)}_{n}-\alpha_{11}(u^{(2)}_{n}-u^{(1)}_{n})-\beta_{11}(u^{(2)}_{n-1}-u^{(1)}_{n})-\alpha_{12}(\phi^{(2)}_{n}-\phi^{(1)}_{n})-\beta_{12}(\phi^{(2)}_{n-1}-\phi^{(1)}_{n})=0
\end{eqnarray}
\begin{eqnarray}
    \label{equ:17}
    m_{n}\Ddot{u}^{(2)}_{n}-\alpha_{11}(u^{(1)}_{n}-u^{(2)}_{n})-\beta_{11}(u^{(1)}_{n+1}-u^{(2)}_{n})-\alpha_{12}(\phi^{(1)}_{n}-\phi^{(2)}_{n})-\beta_{12}(\phi^{(1)}_{n+1}-\phi^{(2)}_{n})=0
\end{eqnarray}
\begin{eqnarray}
    \label{equ:18}
    j_{n}\Ddot{\phi}^{(1)}_{n}-\alpha_{21}(u^{(2)}_{n}-u^{(1)}_{n})-\beta_{21}(u^{(2)}_{n-1}-u^{(1)}_{n})-\alpha_{22}(\phi^{(2)}_{n}-\phi^{(1)}_{n})-\beta_{22}(\phi^{(2)}_{n-1}-\phi^{(1)}_{n})=0
\end{eqnarray}
\begin{eqnarray}
    \label{equ:19}
    j_{n}\Ddot{\phi}^{(2)}_{n}-\alpha_{21}(u^{(1)}_{n}-u^{(2)}_{n})-\beta_{21}(u^{(1)}_{n+1}-u^{(2)}_{n})-\alpha_{22}(\phi^{(1)}_{n}-\phi^{(2)}_{n})-\beta_{22}(\phi^{(1)}_{n+1}-\phi^{(2)}_{n})=0
\end{eqnarray}
\end{subequations}
\end{widetext}
where the coefficients can be derived from the second derivative of the potential energy of each unit cell as detailed in the reference~\cite{miyazawa2022topological}. With the design parameters and mechanical parameters in TABLE~\ref{table:1} and TABLE~\ref{table:2}, the coefficients can be determined to be $\alpha_{11}=\beta_{11}=26850~\mathrm{N\cdot m^{-1}}$, $\alpha_{12}=-\beta_{12}=-819.8~\mathrm{N\cdot rad^{-1}}$ and $\alpha_{22}=\beta_{22}=26.07~\mathrm{N\cdot m \cdot rad^{-1}}$. The periodic boundary condition (Bloch's theorem) is then applied to the linearized equation of motion, such that $u^{(2)}_{n-1}=u^{(2)}_{n}e^{-2ikh_{0}}$, $\phi^{(2)}_{n-1}=\phi^{(2)}_{n}e^{-2ikh_{0}}$, $u^{(1)}_{n+1}=u^{(1)}_{n}e^{2ikh_{0}}$ and $\phi^{(1)}_{n+1}=\phi^{(1)}_{n}e^{2ikh_{0}}$. The band structure $\omega(k)$ can be obtained by solving the eigenvalue equation as a function of the Bloch wave vector $k$ in the first Brillouin zone:
\begin{eqnarray}
    \label{equ:20}
    [\hat{D}(k)+\omega^{2}]U=0
\end{eqnarray}
Here, $\omega$ denotes the angular frequency. $\hat{D}(k)$ is the dynamical matrix as a function of $k$ and $U$ is the corresponding eigenmode $U=[\sqrt{m}u^{(1)},\sqrt{m}u^{(2)},\sqrt{j}\phi^{(1)},\sqrt{j}\phi^{(2)}]^{T}$.

\section{Error analysis in different training ratios}
\label{Appendix B}
\begin{figure}[h]
    \includegraphics[width=0.4\textwidth]{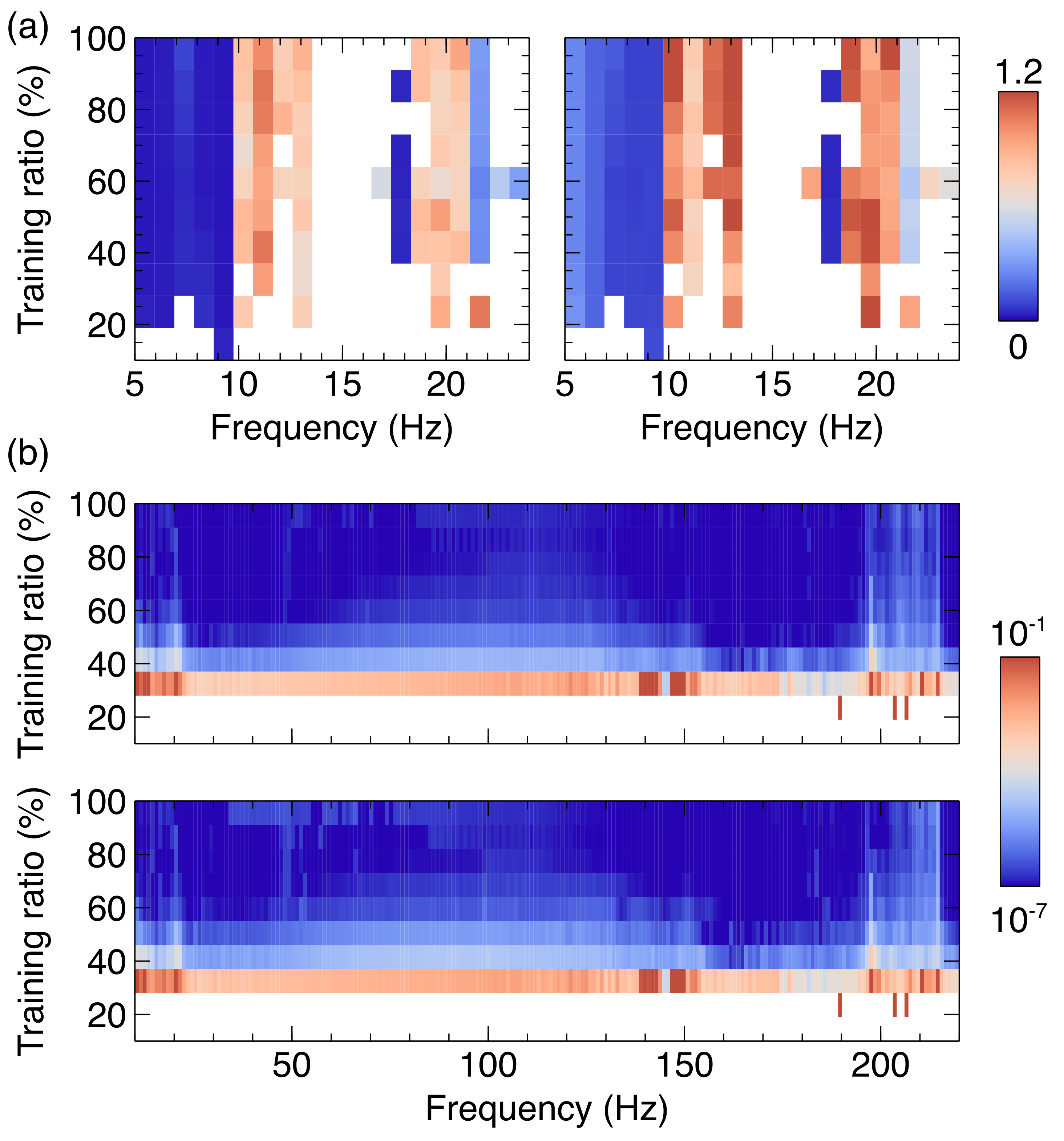}
    \caption{The prediction accuracy for the dual origami structure and the origami chain using different training ratios.
    (a) The relative errors as a function of training ratio of axial and rotational displacement in different frequencies for dual origami structure are shown in the left and right panel, respectively.
    (b) The relative errors as a function of training ratio of axial and rotational displacement in different frequencies for the origami chain are shown in the left and right panel, respectively.
    The white regions indicate the failure of the prediction in (a) and (b).
    \label{fig:fig6}}
\end{figure}

Although the training ratio for the giDMD model of our aforementioned examples is set to be $60\%$, we also investigate the sufficiency of the training ratio to build an effective and reliable model to describe the origami dynamics. The relative error in different training ratios are calculated for two examples (FIG.~\ref{fig:fig6}(a) for the dual origami structure and FIG.~\ref{fig:fig6}(b) for the origami chain). FIG.~\ref{fig:fig6}(a) shows the relative error of axial and rotational displacement in the left and right panels, respectively. It can be seen that $20\%$ training ratio is enough to build the effective giDMD model for most cases of intrawell periodic motion in the dual origami structure. However, giDMD fails to model the interwell periodic motion no matter what the training ratio is. The effective modeling for the chaotic motion can be achieved as early as when the training ratio is $20\%$. The corresponding model cannot precisely describe the intermittent change of chaotic motion as illustrated in FIG.~\ref{fig:fig4}(f) in the main text, resulting in the large relative error. However, it can show the spectral features of the chaotic motion (large low-frequency components) as shown in FIG.~\ref{fig:fig5}(b). Note that large training ratio does not always guarantee the effective giDMD model as shown by the white areas in the region of large training ratios.
FIG.~\ref{fig:fig6}(b) shows the relative error of axial and rotational displacement in the top and bottom panels, respectively. It can be seen that the effective model can be constructed as early as when the training ratio is $30\%$. Besides, the model across frequencies can be built with high accuracy.
\begin{figure}[h!]
    \centering
    \includegraphics[width=0.45\textwidth]{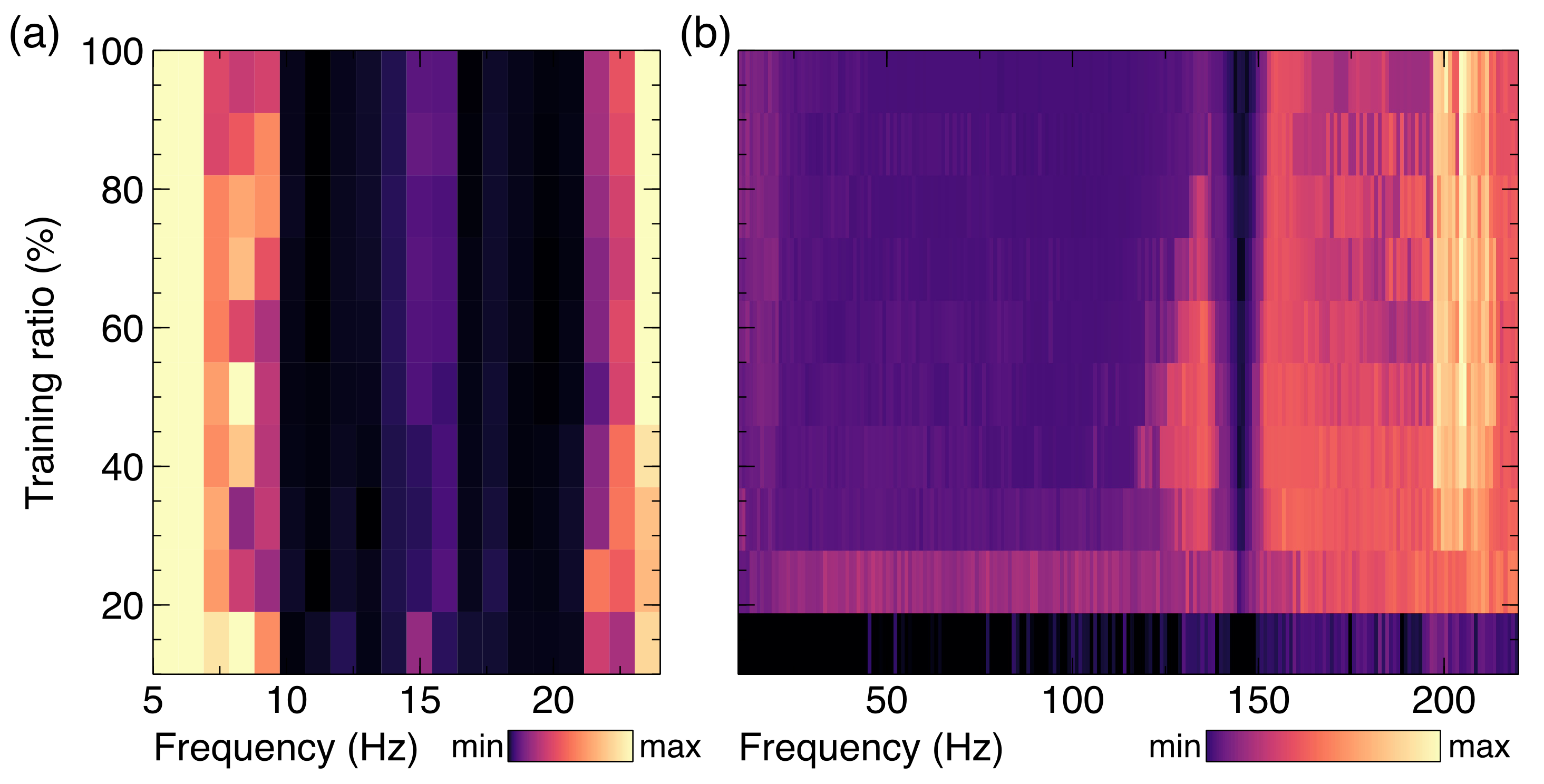}
    \caption{The singular values of $\bm{K}$ matrix for dual origami structure and origami chain using different training ratios.
    (a) The first singular value of $\bm{K}$ matrix as a function of training ratio in different frequencies for dual origami structure.
    (b) The fifth singular value of $\bm{K}$ matrix as a function of training ratio in different frequencies for origami chain.
    \label{fig:fig7}}
\end{figure}
\section{Identification of state in different training ratios}
\label{Appendix C}
As stated in the main text, the `state' of the origami structure may be induced by the singular value of $\bm{K}$ matrix. We further calculate the singular values of $\bm{K}$ matrix in different training ratios to identify the `state' in origami structures. The second singular value of $\bm{K}$ matrix in different training ratios is calculated for the dual origami structure (FIG.~\ref{fig:fig7}(a)), and the fifth singular value is calculated for the origami chain (FIG.~\ref{fig:fig7}(b)). FIG.~\ref{fig:fig7}(a) reveals that the intrawell periodic motion, interwell periodic motion and chaotic motion can be initially identified by distinct values when the training ratio is only $10\%$, although it should be noticed that at this training ratio the giDMD models for different frequencies are not adequate to describe the origami dynamics. In the origami chain, the fifth singular value is used to conduct identification to have a stronger contrast although other singular values also work. The topological boundary state can be identified when the training ratio is $20\%$ even though the model is not effectively built. The study on the training ratio indicates that our approach of giDMD can identify the `state' of origami structures at the early stage, which helps with the early diagnosis of the origami structure.

\bibliography{references.bib}

\end{document}